\def\timestamp{%
Time-stamp: <elem-and-dim.tex: Tuesday 08-07-2003 at 13:20:55 (cest)>}
\def\stripname Time-stamp: <#1 #2>{#2}
\edef\filedate{\expandafter\stripname\timestamp}

\documentclass{amsart}
\newcommand{\Dg}{\operatorname{Dg}}
\newcommand{\ind}{\operatorname{ind}}
\newcommand{\Ind}{\operatorname{Ind}}
\newcommand{\conn}{\operatorname{conn}}
\newcommand{\cut}{\operatorname{cut}}
\newcommand{\partn}{\operatorname{part}}

\let\meet\sqcap
\let\join\sqcup
\newcommand\0{\mathbf{0}}
\newcommand\1{\mathbf{1}}

\let\implies\rightarrow

\renewcommand\newsymbol[5]{%
\DeclareMathSymbol#1{#3}{\ifcase #2\or AMSa\or AMSb\fi}{"#4#5}}
\newsymbol\le      1236         
\newsymbol\ge      123E         
\newsymbol\notle       230A

\theoremstyle{plain}
\newtheorem{proposition}{Proposition}[section]

\begin{document}
\title{Elementarity and dimensions}
\author{Klaas Pieter Hart}
\address{Faculty of Electrical Engineering, Mathematics, and
         Computer Science\\
         TU Delft\\
         Postbus 5031\\
         2600~GA {} Delft\\
         the Netherlands}
\email{K.P.Hart@EWI.TUDelft.NL}
\urladdr{http://aw.twi.tudelft.nl/\~{}hart}
\date{\filedate}


\begin{abstract}
We give an alternative proof of Fedorchuk's recent result that
$\dim X\le\Dg X$ for compact Hausdorff spaces $X$.
We use the L\"owenheim-Skolem theorem to reduce the problem to the
metric case.
\end{abstract}

\maketitle

\section{Introduction}

From the various \emph{topological} notions of dimension that have been
proposed the best-known and most widely used are
$\ind$ and $\Ind$, the small and large inductive dimension, and
$\dim$, the covering dimension.
These capture in various ways the intuition behind dimension.
The inductive dimensions formalize the idea that ``a line is separated
by points, a surface by lines and space by surfaces'',
whereas $\dim$ captures dimension as ``number of directions'',
especially through the theorem on partitions
\cite[7.2.15]{Engelking89}.
These functions assume the same values for all separable metrizable
spaces and assign the correct dimension to Euclidean $n$-space.

In \cite{Brouwer1913} Brouwer proposed another notion of dimension,
Dimensionsgrad ($\Dg$), based on cuts.
It was established only recently in \cite{FedorchukLevinShchepin1999}
that $\Dg$ coincides with the familiar dimension functions on the
class of (locally) compact metric spaces.
Outside of this class $\Dg$ and the dimension functions diverge:
there is, for each~$n$, 
a locally connected complete separable metric space~$X_n$
with $\Dg X_n=1$ and $\dim X_n=n$, see~\cite{FedorchukvanMill2000}.

Recently Fedorchuk proved that $\dim X\le\Dg X$ for compact Hausdorff 
spaces~$X$.
The purpose of this note is to reprove this and Vedenissof's inequality
$\dim X\le\Ind X$ (for normal spaces) by model-theoretic means.

The arguments in this paper seem to indicate that $\Dg$ is somewhat
more complex than the common dimension functions, which may help
to explain why Fedorchuk's proof of his inequality is so much more
involved than the fairly straightforward proof of Vedenissof's
inequality.

\section{Preliminaries}

\subsection{Dimensions}

We repeat the definitions of covering dimension and large inductive dimension.
We say that covering dimension of a normal space~$X$ is at most~$n$,
in symbols $\dim X\le n$, if every finite open cover has a refinement
of order at most $n+1$ (i.e., no point is in more than $n+1$ members
of the refinement).
As usual $\dim X$ is defined to be the minimum~$n$ for which this holds
(or $\infty$ if there is no such~$n$).

The large inductive dimension is defined by recursion: $\Ind X\le n$ means
that between every two disjoint closed sets $A$ and $B$ there is a 
\emph{partition}~$C$ with $\Ind C\le n-1$, where $C$~is a partition between~$A$
and~$B$ if $X\setminus C$ can be written as the union of two disjoint open
sets~$U$ and $V$ with $A\subseteq U$ and $B\subseteq V$.
This recursion starts with $\Ind X\le -1$ if{}f $X=\emptyset$.

The Dimensionsgrad is defined similarly but now $C$ should be a \emph{cut} 
between~$A$ and~$B$, which means that it is closed and meets every continuum 
that intersects both~$A$ and~$B$.

\subsection{Lattices}

In \cite{Wallman38} Wallman showed that to every distributive lattice~$L$
with $\0$ and~$\1$
one can associate a compact $T_1$-space $wL$, its 
\emph{Wallman representation}, with a base for the closed sets that is 
a homomorphic image of~$L$.
The underlying set of~$wL$ is the set of all ultrafilters on~$L$ and for 
every element~$a$ of~$L$ the set $\bar a=\{u\in wL:a\in u\}$ is a basic
closed set in~$wL$.
The homomorphism $a\mapsto\bar a$ is one-to-one if and only if $L$~is
\emph{separative}, 
which means that whenever $a\notle b$ there is $c\le a$ with
$c>\0$ and $c\meet b=\0$.
The space $wL$ is Hausdorff if and only if $L$~is \emph{normal}, 
which means that 
whenever $a\meet b=\0$ there are $f$ and $g$ with $a\meet f=\0$,
$b\meet g=\0$ and $f\join g=\1$.

Unlike the Stone representation for Boolean algebras the Wallman 
representation is not one-to-one.
Certainly every compact $T_1$-space~$X$ is the representation of its own
lattice of closed sets, which we denote by~$2^X$, but one also has 
$X=w\mathcal{B}$ whenever $\mathcal{B}$~is a base for the closed sets
of~$X$ that is closed under finite unions and intersections.
Thus, e.g., the unit interval is also the Wallman representation of
the family of finite unions of closed intervals with rational end points.

\subsection{Elementary sublattices}

Our proofs of Fedorchuk's and Vedenissof's inequalities involve
elementary sublattices of~$2^X$.
A sublattice~$L$ of~$2^X$ is an \emph{elementary} sublattice if every
equation with parameters from~$L$ that has a solution in~$2^X$ already
has a solution in~$L$.
Here `equation' should be taken in a very wide sense.
What we demand is: whenever 
$\phi(x,y,\ldots,a,b,\ldots)$ is a lattice-theoretic formula
with its free variables among $x$, $y$, \dots{} and
its parameters $a$, $b$, \dots{} from~$L$ and if there are $x$, $y$, \dots{}
in~$2^X$ such that $\phi$ holds in~$2^X$ then there are such $x$, $y$, \dots{}
in~$L$.

For example an elementary sublattice of~$2^X$ is automatically separative:
if $a\notle b$ in~$L$ then $(x\le a)\land(x>\0)\land(x\meet b=\0)$
is an equation with parameters --- $a$, $b$ and $\0$ --- from~$L$
and with a solution in~$2^X$, hence there must be a $c\in L$ with
$(c\le a)\land(c>\0)\land(c\meet b=\0)$.

Likewise $L$ must be normal: if $a,b\in L$ and $a\meet b=\0$ then the
equation  
$(a\meet x)\land(b\meet y)\land(x\join y=\1)$
has a solution in~$2^X$, hence there are $f$ and $g$ in~$L$ with
$(a\meet f)\land(b\meet g)\land(f\join g=\1)$.

Below, when proving Fedorchuk's inequality we shall see more complicated
equations/formulas that will involve quantifiers and this is where
the strength of the notion of elementarity will become apparent.

An important result is the L\"owenheim-Skolem
theorem, which says, in our context, that given a subfamily $\mathcal{F}$
of $2^X$ one can always find an elementary sublattice~$L$ of $2^X$
with $\mathcal{F}\subseteq L$ and $|L|\le|\mathcal{F}|\cdot\aleph_0$.
This provides an inroad to a strong version of Marde\v{s}i\'c's
Factorization theorem, see \cite[Theorem 5.3]{HartvanMillPol2002}
for an example of its use and the thesis \cite{vanderSteeg2003} 
for a systematic study of the properties that the factorizing space
inherits from the domain.
A proof of the L\"owenheim-Skolem theorem can be found in
\cite[Section 3.1]{Hodges1997}.

\section{Formulas for dimensions}

\subsection{Covering dimension}

We use Hemmingsen's characterization  
from \cite{Hemmingsen1946} (see also~\cite[Corollary 7.2.14]{Engelking89}) to
make a lattice-theoretic formula that characterizes covering dimension in terms
of closed sets.
The formula, abbreviated $\delta_n$, is
\begin{multline}
(\forall x_1)(\forall x_2)\cdots(\forall x_{n+2})
(\exists y_1)(\exists y_2)\cdots(\exists y_{n+2})\\
\bigl[(x_1\meet x_2\meet\cdots\meet x_{n+2}=\0)\implies
\bigl((x_1\le y_1)\land(x_2\le y_2)\land\cdots\land(x_{n+2}\le y_{n+2})\\
 \land(y_1\meet y_2\meet\cdots\meet y_{n+2}=\0)
 \land(y_1\join y_2\join\cdots\join y_{n+2}=\1)\bigr)\bigr].
\end{multline}
Hemmingsen's theorem simply says that, for compact spaces,
$\dim X\le n$ if and only if the lattice~$2^X$ satisfies 
the formula~$\delta_n$.

A standard shrinking-and-expanding argument will show that for a compact 
Hausdorff space~$X$ one has $\dim X\le n$ if and only if some (every) 
lattice base for its closed sets satisfies the formula~$\delta_n$.

\subsection{Large inductive dimension}\label{subsec.Ind}

The definition of large inductive dimension can be couched in terms of closed
sets quite easily.
A partition~$C$ between two disjoint closed sets~$A$ and~$B$ can be described
by two closed sets~$F$ and~$G$ such that 
$F\cup G=X$, $F\cap A=\emptyset$ and $G\cap B=\emptyset$: the intersection
$F\cap G$~is a partition between~$A$ and~$B$; thus the following
formula $\partn(u,x,y,a)$ states that $u$~is a partition between $x$ and~$y$
in the (sub)space~$a$:
$$
(\exists f)(\exists g)
 \bigl((x\meet f=\0)\land(y\meet g=\0)\land(f\join g=a)\land(f\meet g=u)\bigr).
$$
This enables us to give a recursive definition of a formula~$I_n(a)$
for the large inductive dimension:
\begin{multline}
(\forall x)(\forall y)(\exists u) \\
\bigl[(\bigl((x\le a)\land(y\le a)\land(x\meet y=\0)\bigr)\implies
 \bigl(\partn(u,x,y,a)\land I_{n-1}(u)\bigr)\bigr];
\end{multline}
the recursion starts with $I_{-1}(a)$ abbreviating $a=\0$.

Thus a compact space~$X$ satisfies $\Ind X\le n$ if and only if $2^X$ satisfies
$I_n(\1)$.
More generally, if $X$ has a lattice base $\mathcal{B}$ for its closed sets 
that satisfies $I_n(\1)$ then $\Ind X\le n$; 
this follows readily by induction, once one realizes that
$\{F\in\mathcal{B}:F\subseteq A\}$ is a lattice base for the closed sets 
in the subspace~$A$, when $A\in\mathcal{B}$.

The converse is not true in that not every lattice base for the closed
sets of a space~$X$ with $\Ind X\le n$ must satisfy~$I_n(\1)$.
A simple example is given by the unit interval $[0,1]$ and the 
lattice base generated by the subbase 
$\bigl\{[0,q]:q$ rational$\bigr\}\cup\bigl\{[p,1]:p$ irrational$\bigr\}$.
This lattice does not satisfy $I_n(\1)$ for any~$n$.

\subsection{Dimensionsgrad}\label{subsec.Dg}

As defined above a cut between two (disjoint) closed sets~$A$ and~$B$ 
is a closed set~$C$ such that every continuum from the ambient space 
that intersects~$A$ and~$B$ also intersects~$C$.
If we let $\conn(a)$ abbreviate 
$$
(\forall x)(\forall y)
 \bigl[\bigl((x\meet y=\0)\land(x\join y=a)\bigr)\implies
      \bigl((x=\0)\lor(x=a)\bigr)\bigr],
$$
i.e, ``$a$ is connected'', and use $\cut(u,x,y,a)$ to denote
$$
(\forall v)
\bigl[\bigl((v\le a)\land \conn(v)\land(v\meet x\neq\0)
 \land(v\meet y\neq\0)\bigr)\implies
 (v\meet u\neq\0)\bigr],
$$
i.e., ``$u$ is a cut between $x$ and $y$ in the (sub)space $a$'',
then we get the following recursive definition of a formula $\Delta_n(a)$
for the Dimensionsgrad:
\begin{multline}
(\forall x)(\forall y)(\exists u)\\
\bigl[\bigl((x\le a)\land (y\le a)\land(x\meet y=\0)\bigr)\implies
\bigl(\cut(u,x,y,a)\land \Delta_{n-1}(u) \bigr) \bigr],
\end{multline}
and, as above, $\Delta_{-1}(a)$ denotes $a=\0$.

As with the large inductive dimension one has $\Dg X\le n$ if and only if
$2^X$ satisfies $\Delta_n(\1)$.
The same example as above shows that it is possible to have $\Dg X=1$ 
while some lattice base for the closed does not satisfy $\Delta_n(\1)$
for any~$n$.
It is however also possible that some lattice base for the closed sets of
a space~$X$ satisfies $\Delta_0(\1)$ while $\Dg X>0$.
An example is provided by the unit interval and the lattice base
generated by 
$\bigl\{[0,q]\cup\{q+2^{-n}:n\in\omega\}:q$ rational$\bigr\}
\cup\bigl\{[p,1]\cup\{p-2^{-n}:n\in\omega\}:p$ irrational$\bigr\}$.
This lattice base satisfies $\Delta_0(\1)$ \emph{vacuously}, as it has no
non-trivial connected elements.

\section{Elementarity}

In this section we fix a compact Hausdorff space~$X$ and an
elementary sublattice~$L$ of the lattice~$2^X$, 
with Wallmany representation~$wL$.

\subsection{$\dim wL=\dim X$}

This is by and large well-known but to keep this note self-contained 
we indicate a proof.
By the remarks in the previous section we know that $\dim wL$~is the minimum
natural number~$n$ for which $L$ satisfies $\delta_n$.
Therefore we have to show that $L$ satisfies~$\delta_n$ if and only $2^X$
satisfies~$\delta_n$.
The straightforward part is sufficiency: if $2^X$ satisfies $\delta_n$ then
so does~$L$: every $n+2$-tuple $(x_1,\ldots,x_{n+2})$ from $L$ determines,
via $\delta_n$, an equation that has a solution $(y_1,\ldots, y_{n+2})$
in~$2^X$ and hence in~$L$.
The converse follows by contraposition: the negation of $\delta_n$ is in itself
an equation with only $\0$ and $\1$ as its parameters and unknowns
$x_1$, \dots $x_{n+2}$; if it has a solution in~$2^X$ then it also has a
solution in~$L$.

\subsection{$\Ind wL\le\Ind X$}

As above one deduces that $L$ satisfies $I_n(\1)$ if and only if $2^X$ 
satisfies $I_n(1)$: both $I_n(1)$ and its negation give rise to equations 
with parameters in~$L$ and solutions in~$2^X$, hence in~$L$.
In \ref{subsec.Ind} we have seen that $\Ind wL\le n$ whenever $L$ 
satisfies~$I_n(1)$; this suffices for $\Ind wL\le\Ind X$.

\subsection{$\Dg wL\le\Dg X$}

As above we find that $L$ satisfies $\Delta_n(1)$ if and only if $2^X$
satisfies $\Delta_n(1)$.
However, in \ref{subsec.Dg} we saw that $\Dg wL\le n$ does not follow
automatically from the fact that $L$ satisfies $\Delta_n(1)$.
This shows that a bit more effort will have to go into the proof;
in fact we shall prove the following proposition by induction on~$n$.

\begin{proposition}
Let $X$ be a compact Hausdorff space with $\Dg X\le n$ and
$L$ an elementary sublattice of $2^X$.
Then $\Dg wL\le n$.
\end{proposition}

\begin{proof}
In this proof an element $A$ of $L$ is on the one hand a closed subset of~$X$
and on the other hand a name for a basic closed set in~$wL$; we write
$A_L$ to denote the latter set.

Let $P$ and $Q$ be closed and disjoint sets in~$wL$.
Because $L$~is a lattice base for the closed sets of $wL$ there are 
disjoint $A,B\in L$ with $P\subseteq A_L$ and $Q\subseteq B_L$.

Now in $X$ there a cut $C$ between $A$ and $B$ with $\Dg C\le n-1$,
by elementarity we can assume $C\in L$.
Indeed, apparently there is in~$2^X$ a solution to the equation 
$\cut(x,A,B,\1)\land \Delta_{n-1}(x)$, 
which has parameters in~$L$, hence such a solution must exist in~$L$.

We must show that the closed set $C_L$ represented by~$C$ in~$wL$ is a
cut between $A_L$ and $B_L$ (hence between $P$ and $Q$) 
and that $\Dg(C_L)\le n-1$.

The latter follows by induction because $C_L$ is the Wallman representation
of the lattice $\{x\in L:x\subseteq C\}$ and because this lattice
is an elementary sublattice of $\{x\in 2^X:x\subseteq C\}$.

To prove the former assume $K$~is a closed set in~$wL$ that meets $A_L$ 
and~$B_L$ but not~$C_L$.
Take $H\in L$ with $K\subseteq H_L$ and $H\cap C=\emptyset$.
Observe that $H$~is not connected because it intersects both~$A$ and~$B$
but not~$C$.
One can therefore apply elementarity to the formula $\lnot\conn(H)$ to 
find non-zero disjoint elements $F$ and $G$ of~$L$ with $H=F\cup G$.
Then $H_L$~is the disjoint union of~$F_L$ and~$G_L$; this does not
help in proving $K$ disconnected however, as it is quite possible 
that $K\subseteq F_L$ or $K\subseteq G_L$.
We shall have to choose $F$ and $G$ with extra care.

We use the fact that, in~$X$, no component of~$H$ meets both~$A$ and~$B$.
Because the decomposition of~$H$ into its components is upper-semicontinuous
\cite[6.2.21]{Engelking89}
it follows that we can find two disjoint closed sets $F$ and $G$ such that
$F\cup G=H$, $A\cap H\subseteq F$ and $B\cap H\subseteq G$. 
Again, elementarity dictates that there are such $F$ and $G$ in~$L$.

Now work in $wL$:
as $K\subseteq H_L$, we know that 
$\emptyset\neq K\cap A_L\subseteq H_L\cap A_L\subseteq F_L$ and
$\emptyset\neq K\cap B_L\subseteq H_L\cap B_L\subseteq G_L$.
But this implies $K$~is not connected, as $K\subseteq H_L=F_L\cup G_L$
and $F_L\cap G_L=\emptyset$.

This shows, by contraposition, that $C_L$ is indeed a cut between $A_L$ 
and~$B_L$.
\end{proof}

\subsection{Proof of $\dim X\le \Dg X$} \label{subsec.dim<=Dg}

By the L\"owenheim-Skolem theorem one can find a \emph{countable} elementary 
sublattice $M$ of $2^X$.
The Wallman representation $wM$ of this lattice is compact and metrizable.

The theorem from \cite{FedorchukLevinShchepin1999} says that
$\dim wM=\Dg wM$.
Combined with the equality $\dim wM=\dim X$ and the inequality $\Dg wM\le\Dg X$
this establishes $\dim X\le\Dg X$.

\section{Concluding remarks}

As every partition between two closed sets is also a cut between these sets
one gets the inequality $\Dg X\le \Ind X$ for normal spaces without any
real effort.
A consequence of Fedorchuk's inequality is Vedenissof's inequality 
$\dim X\le \Ind X$ for compact spaces
\cite[7.2.8]{Engelking89}.
The L\"owenheim-Skolem method can also be used to prove this directly:
with the notation as in~\ref{subsec.dim<=Dg} one has
$\dim X=\dim wM=\Ind wM\le\Ind X$ (an application of the \v{C}ech-Stone
compactification allows one to extend this to all normal spaces).

As remarked in the introduction the standard proof of $\dim X\le\Ind X$ 
is fairly straightforward, whereas Fedorchuk's proof in \cite{Fedorchuk2003} 
of $\dim X\le\Dg X$ is longer and needs a closing-off argument to find 
a good cut.
This difference is also apparent in the proofs in the present paper:
in both cases the first step was to produce a (countable) lattice $M$ 
that satisfies $I_n(\1)$ or $\Delta_n(\1)$ respectively.
The second step was to deduce that $\Ind wM\le n$ or $\Dg wM\le n$ 
respectively.
In either case the formula produced a candidate partition or cut;
the problem then was to show that this set was indeed a partition or cut
in the space~$wM$.
This is easy in the case of a partition: once the closed sets $F$ and $G$
are found we are done.
In the case of a cut we only know that our set meets the connected elements
of the base~$M$ that meet~$A$ and~$B$; we need to know that the same holds
for all continua in~$wM$.
This is where elementarity was used once more: it ensured that $M$~already
contained enough connected elements for the proof to go through.
The reason for the perceived unwieldiness of $\Dg$ therefore seems to stem
from the hidden universal quantifiers in the formula $\cut(u,x,y,a)$


\begin{thebibliography}{10}

\bibitem{Brouwer1913}
L.~E.~J. Brouwer, \emph{{\"{U}}ber den nat\"{u}rlichen {D}imensionsbegriff},
  Journal f\"ur die reine und angewandte Mathematik \textbf{142} (1913),
  146--152.

\bibitem{Engelking89}
Ryszard Engelking, \emph{General {Topology}. {Revised} and completed edition},
  Sigma Series in Pure Mathematics, no.~6, Heldermann Verlag, Berlin, 1989.
  \MR{91c:54001}

\bibitem{Fedorchuk2003}
V.~V. Fedorchuk, \emph{On the {Brouwer} dimension of compact spaces},
  Mathematical Notes \textbf{73} (2003), 271--279, Russian original:
  Matematicheskie Zametki \textbf{73} (2003) 295--304.

\bibitem{FedorchukLevinShchepin1999}
V.~V. Fedorchuk, M.~Levin, and E.~V. Shchepin, \emph{On the {B}rouwer
  definition of dimension}, Uspekhi Matematicheskikh Nauk \textbf{54} (1999),
  no.~2(326), 193--194, translation in Russian Mathematical Surveys \textbf{54}
  (1999), no. 2, 432--433. \MR{2000g:54068}

\bibitem{FedorchukvanMill2000}
V.~V. Fedorchuk and J.~van Mill, \emph{Dimensionsgrad for locally connected
  {Polish} spaces}, Fundamenta Mathematicae \textbf{163} (2000), 77--82.
  \MR{2001a:54045}

\bibitem{HartvanMillPol2002}
Klaas~Pieter Hart, Jan van Mill, and Roman Pol, \emph{Remarks on hereditarily
  indecomposable continua}, Proceedings of the 15th Summer Conference on
  General Topology and its Applications/1st Turkish International Conference on
  Topology and its Applications (Oxford, OH/Istanbul, 2000), vol.~25, 2000,
  pp.~179--206 (2002). \MR{1 925 683}

\bibitem{Hemmingsen1946}
Erik Hemmingsen, \emph{Some theorems in dimension theory for normal {H}ausdorff
  spaces}, Duke Mathematical Journal \textbf{13} (1946), 495--504. \MR{8,334e}

\bibitem{Hodges1997}
Wilfrid Hodges, \emph{A shorter model theory}, Cambridge University Press,
  Cambridge, 1997. \MR{98i:03041}

\bibitem{vanderSteeg2003}
Berd van~der Steeg, \emph{Models in topology}, Ph.D. thesis, TU Delft, 2003.

\bibitem{Wallman38}
Henry Wallman, \emph{Lattices and topological spaces}, Annals of Mathematics
  \textbf{39} (1938), 112--126.

\end{thebibliography}

\providecommand{\bysame}{\leavevmode\hbox to3em{\hrulefill}\thinspace}
\providecommand{\MR}{\relax\ifhmode\unskip\space\fi MR }
\providecommand{\MRhref}[2]{%
  \href{http://www.ams.org/mathscinet-getitem?mr=#1}{#2}
}
\providecommand{\href}[2]{#2}

\end{document}